\documentclass[10pt]{article}
\usepackage{authblk}
\usepackage{graphicx}
\usepackage{amsmath}

\begin{document}

\title{Lie Symmetry Analysis of the Black-Scholes-Merton Model for European
Options with Stochastic Volatility}
\author[1]{A Paliathanasis\thanks{%
paliathanasis@na.infn.it}}
\author[2]{K Krishnakumar\thanks{%
krishapril09@gmail.com}}
\author[2]{KM Tamizhmani\thanks{%
kmtmani54@gmail.com}}
\author[3,4,5]{ PGL Leach\thanks{%
leach@ucy.ac.cy}}

\affil[1]{Instituto de Ciencias F\'{\i}sicas y Matem\'{a}ticas, Universidad Austral de Chile, Valdivia, Chile}
\affil[2]{Department of Mathematics, Pondicherry University, Kalapet Puducherry 605 014, India}
\affil[3]{Department of Mathematics and
Institute of Systems Science, Durban
University of Technology, PO Box 1334, Durban 4000, Republic of South Africa}
\affil[4]{School of Mathematics, Statistics and Computer Science, University
of KwaZulu-Natal, Private Bag X54001, Durban 4000, Republic of South Africa}
\affil[5]{Department of Mathematics and Statistics, University of Cyprus,
Lefkosia 1678, Cyprus}

\renewcommand\Authands{ and }
\maketitle

\begin{abstract}
We perform a classification of the Lie point symmetries for the
Black--Scholes--Merton Model for European options with stochastic volatility, $\sigma$, in which the last is defined by a stochastic differential equation
with an Orstein--Uhlenbeck term. In this model, the value of the option is
given by a linear (1 + 2) evolution partial differential equation in which
the price of the option depends upon two independent variables, the value of
the underlying asset, $S$, and a new variable, $y$. We find that for arbitrary functional form of the
volatility, $\sigma(y)$, the (1 + 2) evolution equation always admits two Lie point
symmetries in addition to the automatic linear symmetry and the infinite number of solution
symmetries. However, when $\sigma(y)=\sigma_{0}$ and as the price of the
option depends upon the second Brownian motion in which the volatility is
defined, the (1 + 2) evolution is not reduced to the Black--Scholes--Merton
Equation, the model admits five Lie point symmetries in addition to the linear symmetry and
the infinite number of solution symmetries.  We apply the
zeroth-order invariants of the Lie symmetries and we reduce the (1 + 2)
evolution equation to a linear second-order ordinary differential equation.
Finally, we study two models of special interest, the Heston model and the
Stein--Stein model.
\end{abstract}


\section{Introduction}

The Black--Scholes--Merton Model for European options is based upon some
Ansatz for the stock price. Specifically, the process for the stock price is
characterized by continuity, and it has the ability to hedge continuously with
transaction costs and has constant volatility \cite{bsch01,bsch02,Merton}. 

In the Black--Scholes--Merton Model, the price of a financial asset is given by the soluton of
the stochastic differential equation%
\begin{equation}
dS_{t}=rS_{t}dt+\sigma S_{t}dW_{t}  \label{bs.00}
\end{equation}%
where $W_{t}$ is a Brownian motion, and the value $u=u\left( t,S\right) $ of
the option is given by the solution of the $(1+1)$ evolution equation,
\begin{equation}
\frac{1}{2}\sigma ^{2}S^{2}u_{,SS}+rSu_{,S}-ru+u_{,t}=0  \label{bs.01}
\end{equation}%
in which $t$ is time, $S$ is the current value of the underlying asset, for
example a stock price, and $r$ is the rate of return on a safe investment. The
value of the option is subject to the satisfaction of the terminal
condition, $u\left( T,S\right) =U\left( S\right) $, when $t=T$. Finally, $%
\sigma $ is the volatility of the model.

The Black--Scholes--Merton Model assumes constant volatility $\sigma $.
However, in real problems, $\sigma $ is not a constant. One possible
generalisation of the model, Equation (\ref{bs.01}), is to consider that the
volatility depends upon the time, $t$, and on the value of the stock, $S$, \textit{i.e.}, $%
\sigma =\sigma \left( t,S\right) $. It has been proposed that $\sigma $ is \mbox{a
function} of a mean Orstein--Uhlenbeck process \cite{Sophocleous 11 a}.

Consider that $\sigma =f\left( y\right) $, where $y$ is given by the
stochastic differential equation with the Orstein--Uhlenbeck term \cite%
{finBook,finBook2,hull}:
\begin{equation}
dy_{t}=\alpha \left( m-y_{t}\right) dt+\beta d\hat{Z}_{t}  \label{sv.03}
\end{equation}%

The new Brownian motion $\hat{Z}_{t}$ can be correlated with $W_{t}$ and be
expressed as follows:
\begin{equation}
\hat{Z}_{t}=\rho W_{t}+\sqrt{1-\rho ^{2}}Z_{t}  \label{sv.04}
\end{equation}%
in which $Z_{t}$ describes a Brownian motion independent of $W_{t}$ and
$\rho$ is the correlation factor with values $\vert \rho  \vert
\leq 1$.

Hence, the Black-Scholes Equation (\ref{bs.01}) in the case of stochastic
volatility is modified and the value $u$ of the option is given by the $(1+2)$
evolution equation%
\begin{equation}
\left( \hat{M}_{1}+\hat{M}_{2}+\hat{M}_{3}+\hat{M}_{4}\right) u\left(
t,S,y\right) =0  \label{sv.05}
\end{equation}%
where the operators, $\hat{M}_{1},$~$\hat{M}_{2},~\hat{M}_{3},~\hat{M}_{4},$
are defined as follows:%
\begin{equation}
\hat{M}_{1}=\frac{1}{2}f^{2}\left( y\right) S^{2}\frac{\partial ^{2}}{%
\partial S^{2}}+rS\frac{\partial }{\partial S}-r+\frac{\partial }{\partial t}%
  \label{sv.06}
\end{equation}%
\begin{equation}
\hat{M}_{2}=\rho \beta Sf\left( y\right) \frac{\partial ^{2}}{\partial
S\partial y}~,~\hat{M}_{3}=-\beta \Lambda \left( t,S,y\right) \frac{\partial
}{\partial y}, \quad \mbox{\rm and}  \label{sv.07}
\end{equation}%
\begin{equation}
\hat{M}_{4}=\frac{1}{2}\beta ^{2}\frac{\partial ^{2}}{\partial y^{2}}+\alpha
\left( m-y\right) \frac{\partial }{\partial y}  \label{sv.08}
\end{equation}%

The function $\Lambda \left( t,S,y\right) $ is%
\begin{equation}
\Lambda \left( t,S,y\right) =\rho \frac{\mu -r}{f\left( y\right) }+\gamma
\left( t,S,y\right) \sqrt{1-\rho ^{2}}  \label{sv.09}
\end{equation}%
and $u\left( t,S,y\right) $ satisfies the terminal condition $u\left(
T,S,y\right) =U\left( S\right) $ at time $t=T$.

The operator $\hat{M}_{1}$ gives the Black--Scholes--Merton Equation (\ref%
{bs.01}) with volatility $\sigma =f\left( y\right) $, $\hat{M}_{2}$
expresses the correlation term between the two Brownian motions, $W_{t}$ and
$\hat{Z}_{t}$, of the European option and of the volatility, respectively,
and $\hat{M}_{4}$ is the Orstein--Uhlenbeck process term. Finally. the term $~%
\hat{M}_{3}$, the so called premium term, expresses the market price of the
volatility risk \cite{finBook2}. The function $\gamma \left( t,S,y\right) $
in Equation (\ref{sv.09}) is the risk-premium factor which drives the volatility and
follows from the second Brownian motion, $Z_{t}$, where in the case of
absolute correlation, \textit{i.e.}, $\left\vert \rho \right\vert =1$, $\gamma
\left( t,S,y\right) $ does not play any role in the model. The first term of
the {\it rhs} side of Equation (\ref{sv.09}) is called the excess return-to-risk \mbox{ratio \cite%
{finBook2}}. \mbox{The statistical} importance of stochastic volatility has been
confirmed in \cite{ander}.

The purpose of this work is the study of the Black--Scholes--Merton Model with
stochastic volatility, Equation (\ref{sv.05}), by using the method of group invariant
transformations, specifically the Lie (point) symmetries of the equation. The importance of
Lie symmetries is that they provide a systematic method to facilitate the
solution of differential equations because they provide first-order
invariants which can be used to reduce the differential equations. Moreover,
Lie symmetries can be used for the classification of differential equations.
Furthermore, we can extract important information for the differential
equation, consequently for the model, from the group of invariant
\mbox{transformations admitted}.

The first application of the Lie symmetries in financial modeling was
performed by Gazizov \& Ibragimov in \cite{ibrag}.  They studied the admitted
group of invariant transformations for the Black--Scholes--Merton Equation (%
\ref{bs.01}), with constant volatility and they proved that Equation (\ref%
{bs.01}) admits as Lie symmetries the elements of Lie algebra,
{$\left\{
A_{3,8}\oplus _{s}A_{3,1}\right\} \oplus _{s}\infty A_{1}$}  (
{In the Mubarakzyanov Classification Scheme} \cite%
{Morozov58a,Mubarakzyanov63a,Mubarakzyanov63b,Mubarakzyanov63c}). This means
that Equation (\ref{bs.01}) is maximally symmetric and according to the Theorem of
Sophus Lie \cite{lie} there exists a transformation on the space of
variables $\left\{ t,S,u\right\} $ in which Equation (\ref{bs.01}) can be written in
the form of the heat equation. The last was an important result because the
mathematical methods from physical science can be used for the study of
differential equations in financial mathematics. A similar result has been
found for the one-factor model of commodities \cite{Leach08}, which means
that the three different equations, the heat equation, the
Black--Scholes--Merton equation and the one-factor model of commodities
equation, are equivalent at the mathematical level even if they describe
different subjects.

In recent years, Lie symmetries have covered a big range of applications in
financial mathematics. For instance, the group invariants of the
Cox--Ingersoll--Ross Pricing Equation have been studied \mbox{in \cite{ohara1}} and the
nonlinear Merton model in \cite{ohara2}. As far as concerns the Asian option,
a Lie symmetry classification has been performed in \cite{asian}. As for
generalisations of the Black--Scholes--Merton Model, the Lie
symmetries and the reduction process of the nonautonomous model can be found
in \cite{soph,Tamizhmani}, while another generalisation of Equation (\ref%
{bs.01}) with a \textquotedblleft source\textquotedblright\ was studied in
\cite{dimas2}.

Furthermore, in \cite{paper1,paper2}, the symmetry analysis of the space- and
time-dependent one-factor model of commodities and of the nonautonomous
two-dimensional Black--Scholes--Merton Equations were performed. For other
applications of Lie symmetries in financial mathematics, see, for instance,
\cite{fin1,fin2,fin3}, and references therein.

The stochastic volatility model, Equation (\ref{sv.05}), is a $(1+2)$ evolution
equation.  Below, we perform a symmetry analysis and we determine the group
invariant solutions. \ In particular, we restrict our analysis to the model
in which the risk premium factor vanishes without necessarily $\left\vert
\rho \right\vert =1$ and from Equation (\ref{sv.09}), only the term which expresses the
return-to-risk ratio survives. Moreover, we study\ \mbox{two models} for European
options with stochastic volatility, the Heston model \cite{heston} and the
Stein--Stein model \cite{steinstein}. The latter is a model without
correlation between the \mbox{two Brownian} motions, $W_{t}$ and $\hat{Z}_{t},~$
\textit{i.e.}, $\rho =0$ in Equation (\ref{sv.04}). The plan of the paper is as follows.

In Section \ref{section2}, we give the basic properties and definitions for
the Lie point symmetries of differential equations and we perform the
symmetry classification for our model. We find that Equation (\ref{sv.05})
without the risk premium factor is always invariant under the $\left\{
3A_{1}\right\} \\oplus _{s}\infty A_{1}$ Lie algebra. However, when $f\left(
y\right) $ is constant, Equation (\ref{sv.05}) is invariant under a larger
Lie algebra. The application of the Lie symmetries to Equation (\ref{sv.05})
can be found in Section \ref{group}, in which we reduce the $\left(
1+2\right) $ evolution equation by using the zeroth-order invariants
provided by the Lie symmetries and we derive invariant solutions. \mbox{In
Sections \ref{Heston}} and \ref{steinstein} we study two models of stochastic
volatility for European options, the Heston model and the Stein--Stein model,
respectively. For these two models, we find that both are invariant under the
Lie algebra $\left\{ 3A_{1}\right\} \oplus _{s}\infty A_{1}$, and we apply
the Lie symmetries to solve the equations of the two models. For the Heston
model, the closed-form solution is expressed in terms of Kummer Functions,
whereas for the Stein--Stein model, the closed-form solution is expressed in
terms of Hypergeometric Functions. Furthermore, we give some numerical
solutions for the two models. Finally, in Section \ref{conclusions}, we
discuss our results and draw our conclusions

\section{Lie Symmetry Analysis}

\label{section2}

We consider the Black--Scholes--Merton Equation with stochastic volatility
governed by the evolution Equation (\ref{sv.05}) for which the premium term
depends only upon the return-to-risk ratio. For a time-independent rate-of-return, Equation (\ref{sv.05}) becomes
\begin{eqnarray}
\mathcal{H} &:&0=\frac{1}{2}f^{2}\left( y\right) S^{2}u_{,SS}+\rho \beta
Sf\left( y\right) u_{,Sy}+\frac{1}{2}\beta ^{2}u_{,yy}+  \\
&&+rSu_{,S}+\left[ \alpha \left( m-y\right) -\beta \rho \frac{\mu -r}{%
f\left( y\right) }\right] u_{,y}-ru+u_{,t}  \label{sv.10}
\end{eqnarray}

Let $\Phi $ be the map of an one-parameter point transformation such as
\begin{equation}
\Phi \left( u\left( t,S,y\right) \right) =u^{\prime }\left( t^{\prime
},S^{\prime },y^{\prime }\right)  \label{sv.10a}
\end{equation}%
$\ $with infinitesimal transformation ({$\varepsilon $ is the
parameter of smallness.})%
\begin{eqnarray}
t^{\prime } &=&t+\varepsilon \xi ^{1}\left( t,S,y,u\right)   \label{sv.11}
\\
S^{\prime } &=&S+\varepsilon \xi ^{2}\left( t,S,y,u\right)   \label{sv.12}
\\
y^{\prime } &=&y+\varepsilon \xi ^{3}\left( t,S,y,u\right)   \label{sv.13}
\\
u^{\prime } &=&u+\varepsilon \eta \left( t,S,y,u\right)  \label{sv.15}
\end{eqnarray}%
and generator
\begin{equation}
X=\frac{\partial t^{\prime }}{\partial \varepsilon }\partial _{t}+\frac{%
\partial S^{\prime }}{\partial \varepsilon }\partial _{S}+\frac{\partial
y^{\prime }}{\partial \varepsilon }\partial _{y}+\frac{\partial u^{\prime }}{%
\partial \varepsilon }\partial _{u}  \label{sv.16}
\end{equation}

Consider now that $u\left( t,S,y\right) $ is a solution of  Equation (\ref{sv.10}) and
under the map $\Phi $, \mbox{Equation (\ref{sv.10a})}, $u^{\prime }\left( t^{\prime
},S^{\prime },y^{\prime }\right) $ is also a solution of Equation (\ref{sv.10}). Then,
we say that the generator $X,$ of the infinitessimal transformation of the
one-parameter point transformation, $\Phi$, \ is a Lie (point) symmetry of (%
\mbox{Equation \ref{sv.10})} and Equation (\ref{sv.10}) is invariant under the action of the
map $\Phi $. That means that there exists a function $\psi $ such that the
following condition holds \cite{Bluman}%
\begin{equation}
X^{\left[ 2\right] }\left( \mathcal{H}\right) =\psi \mathcal{H},~mod\left(
\mathcal{H}\right) =0  \label{sv.17}
\end{equation}%
or, equivalently,%
\begin{equation}
X^{\left[ 2\right] }\left( \mathcal{H}\right) =0  \label{sv.17a}
\end{equation}%
where $X^{\left[ 2\right] }$ is the second prologation/extension of $X$ in
the space of variables $\left\{
t,S,y,u,u_{,S},u_{,y},u_{,SS},u_{,Sy},u_{yy}\right\} $. \ Specifically $X^{%
\left[ 2\right] }$ is defined from the following formula
\begin{equation}
X^{[2]}=X+\eta _{i}^{A}\partial _{u^{A}}+\eta _{ij}^{A}\partial _{u_{ij}^{A}}
\end{equation}%
where $\eta _{i}^{A},~\eta _{ij}^{A}$ are given by the relations%
\begin{equation}
\eta _{i}^{A}=\eta _{,i}^{A}+u_{,i}^{B}\eta _{,B}^{A}-\xi
_{,i}^{j}u_{,j}^{A}-u_{,i}^{A}u_{,j}^{B}\xi _{,B}^{j}~
\end{equation}%
and
\begin{align}
\eta _{ij}^{A}& =\eta _{,ij}^{A}+2\eta _{,B(i}^{A}u_{,j)}^{B}-\xi
_{,ij}^{k}u_{,k}^{A}+\eta _{,BC}^{A}u_{,i}^{B}u_{,j}^{C}-2\xi
_{,(i|B|}^{k}u_{j)}^{B}u_{,k}^{A}  \notag \\
& -\xi _{,BC}^{k}u_{,i}^{B}u_{,j}^{C}u_{,k}^{A}+\eta
_{,B}^{A}u_{,ij}^{B}-2\xi _{,(j}^{k}u_{,i)k}^{A}-\xi _{,B}^{k}\left(
u_{,k}^{A}u_{,ij}^{B}+2u_{(,j}^{B}u_{,i)k}^{A}\right)
\end{align}

The importance of the existence of a Lie symmetry for a partial differential
equation is that from the associated Lagrange's system,%
\begin{equation}
\frac{dx^{i}}{\xi ^{i}}=\frac{du^{A}}{\eta ^{A}}
\end{equation}%
Zeroth-order invariants,~$U^{\left[ 0\right] }\left( x^{k},u^{A}\right) $, can
be determined which can be used to reduce the number of the
independent variables of the differential equation.

In the following, we perform a classification of the Lie symmetries of
Equation (\ref{sv.10}). \mbox{Function $f\left( y\right) $} is defined by the
requirement that Equation (\ref{sv.10}) admit Lie symmetries. The latter
requirement can be seen as a geometric selection rule as the Lie symmetries
are generated from the elements of the Homothetic Algebra \cite{anJGP} of
the (pseudo)Riemannian space, which defines the Laplace operator in the $(1+2)$
evolution Equation (\ref{sv.10}). In our case, the (pseudo)Riemannian
manifold is defined by the Brownian motions, $W_{t},~\hat{Z}_{t},$ of the
stock price, $S$, and of the volatility, $\sigma ~$, respectively.

Before we proceed with the symmetry analysis, we remark that Equation (\ref{sv.10}) is
a linear equation which means that it always admits the linear symmetry, $%
X_{u}=u\partial _{u}$ and the infinite-dinensional abelian subalgebra of solutions, $X_{b}=b\left(
t,S,y\right) \partial _{u}$, where function $b\left( t,S,y\right) ,$ is a
solution of the original Equation (\ref{sv.10}) \cite{bluman2}.

\subsection{Classification}

From the symmetry condition Equation (\ref{sv.17}), we get a system of thirty-one
equations (
{for the derivation of the system, we used the symbolic package SYM of
Mathematica} \cite{Dimas05a,Dimas06a})
in which the solution of the system
gives the form of the generator Equation (\ref{sv.16}) of the transformation \mbox{Equation (\ref
{sv.10a})}, that transforms solutions into solutions.  From the latter
system, we have the following results.

For arbitrary function, $f\left( y\right) $, Equation (\ref{sv.10}) admits
the Lie symmetries
\begin{equation}
X_{1}=\partial _{t}~,~X_{2}=S\partial _{S}  \label{sv.18}
\end{equation}%
plus the vector fields $X_{u},~X_{b}$. The algebra in which the Lie
symmetries form is the$\left\{ 3A_{1}\right\} \oplus _{s}\infty A_{1}$.

When $f\left( y\right) =f_{0}$, Equation (\ref{sv.10}) admits
the Lie symmetries
\begin{equation}
\bar{X}_{1}=\partial _{t}~,~\bar{X}_{2}=S\partial _{S}~,~\bar{X}%
_{3}=e^{-\alpha t}\partial _{y}~  \label{sv.19}
\end{equation}%
\begin{equation}
\bar{X}_{4}=2f_{0}tS\partial _{S}+2\frac{\beta }{\alpha }f_{0}\rho \partial
_{y}+\left( \left( f_{0}^{2}-2r\right) t+2\ln S\right) u\partial _{u}
\label{sv.20}
\end{equation}%
\begin{equation}
\bar{X}_{5}=e^{\alpha t}\left( 2f_{0}^{2}\beta \rho S\partial _{S}+\beta
^{2}f_{0}\partial _{y}\right) +\left( 2\alpha f_{0}\left( y-m\right) +2\beta
\left( \mu -r\right) \right) u\partial _{u}  \label{sv.22}
\end{equation}%
plus the vector fields $X_{u},~X_{b}$. The Lie Brackets of the Lie algebra
are given in Table \ref{table1}.

\begin{table}[tpb] \centering%
\caption{Lie brackets of the Lie symmetries of Equation (\ref{sv.10}) for
$f(y)=f_{0}$.}%
\begin{tabular}{ccccccc}
$\left[ \mathbf{\bar{X}}_{I},\mathbf{\bar{X}}_{J}\right] $ & $\mathbf{\bar{X}%
}_{1}$ & $\mathbf{\bar{X}}_{2}$ & $\mathbf{\bar{X}}_{3}$ & $\mathbf{\bar{X}}%
_{4}$ & $\mathbf{\bar{X}}_{5}$ & $\mathbf{X}_{u}$ \\ 
$\mathbf{\bar{X}}_{1}$ & $0$ & $0$ & $-a\bar{X}_{3}$ & $2f_{0}^{2}\bar{X}%
_{2}+\left( f_{0}^{2}-2r\right) X_{u}$ & $\alpha \bar{X}_{5}$ & $0$ \\
$\mathbf{\bar{X}}_{2}$ & $0$ & $0$ & $0$ & $2X_{u}$ & $0$ & $0$ \\
$\mathbf{\bar{X}}_{3}$ & $\alpha \bar{X}_{3}$ & $0$ & $0$ & $0$ & $2\alpha
f_{0}X_{u}$ & $0$ \\
$\mathbf{\bar{X}}_{4}$ & $-2f_{0}^{2}\bar{X}_{2}-\left( f_{0}^{2}-2r\right)
X_{u}$ & $2X_{u}$ & $0$ & $0$ & $0$ & $0$ \\
$\mathbf{\bar{X}}_{5}$ & $-a\bar{X}_{5}$ & $0$ & $-2\alpha f_{0}X_{u}$ & $0$
& $0$ & $0$ \\
$\mathbf{X}_{u}$ & $0$ & $0$ & $0$ & $0$ & $0$ & $0$ \\
\end{tabular}%
\label{table1}%
\end{table}%

We remark that the two-factor model of commodities is invariant under the
same algebra of point transformations \cite{Leach08,paper2}. That is an
expected result because the two-factor model of commodities follows from the
one-factor model in which the second factor, product, follows an
Orstein--Uhlenbeck process. Moreover, as we discussed in the Introduction,
the one-factor model is maximally symmetric just like the Black--Scholes--Merton
Equation.

On the other hand $f\left( y\right) =f_{0}$ means that the volatility $%
\sigma $ is constant. However, the second Brownian motion, $\hat{Z}_{t},$ in
the space wherein $\sigma $ is defined, interacts with the Brownian motion $%
W_{t}$ and modifies the Black--Scholes--Merton Model. However, in the
case for which the correlation $\rho $ vanishes, \textit{i.e.}, $\rho =0$,
Equation (\ref{sv.10}) is not reduced to Equation (\ref{bs.01}) but only when the
Orstein--Uhlenbeck process is identically zero, that is, $\beta =0,~\alpha =0$%
. Otherwise, the price $u$ depends upon the \mbox{Orstein--Uhlenbeck process}.

We continue with the reduction of Equation (\ref{sv.10}) by applying the
zeroth-order Lie invariants. Furthermore for every reduced equation we study
the Lie symmetries.

\section{Group Invariant Solutions}

\label{group}

In this section, we apply the Lie symmetries in order to reduce Equation (\ref%
{sv.10}). We study the \mbox{two cases}, $f\left( y\right) =f_{0},$ and $f\left(
y\right) $ to be an arbitrary function. In order to perform the reduction
and the later equation to give a solution of the original problem, there
should be a constraint between the Lie symmetry vector and the terminal
condition. \ However, we perform the reduction without considering the
terminal condition at the moment because the initial conditions can be
modified from different options. As far as the invariant solutions of the
Black--Scholes--Merton Equation (\ref{bs.01}) are concerned, \mbox{see \cite{inv1}}.

\subsection{Arbitrary Function $f\left( y\right) $}

\label{arbi}

For an arbitrary functional form of $f\left( y\right) $, as we saw above,
Equation (\ref{sv.10}) admits three Lie point symmetries in addition to the
infinite number of solution symmetries. The last cannot be used for the
reduction. Hence, we do not consider them. \ Moreover, a solution in which $u$
does not depend upon one of the independent variables is not an acceptable
solution, that is, the static solution is of no interest. Therefore, we
perform reductions with the symmetry vectors $Y_{1}=X_{1}+\kappa _{1}X_{u}$,$%
~Y_{2}=X_{1}+\kappa_{2}X_{u}$ and $Y_{12}=X_{1}+cX_{2}+\kappa _{3}X_{u}$.

Reduction with respect to the Lie invariants of the symmetry vector $Y_{1}$
gives
\begin{equation}
u\left( t,S,y\right) = \exp[\kappa _{1}t]v\left( S,y\right)  \label{sv.23a}
\end{equation}%
where $v\left( S,y\right) $ satisfies the equation%
\begin{eqnarray}
0 &=&\frac{1}{2}f^{2}\left( y\right) S^{2}v_{,SS}+\rho \beta Sf\left(
y\right) v_{,Sy}+\frac{1}{2}\beta ^{2}v_{,yy}  \notag \\
&&+rSv_{,S}+\left[ \alpha \left( m-y\right) -\beta \rho \frac{\mu -r}{%
f\left( y\right) }\right] v_{,y}-\left( r-\kappa _{1}\right) v
\label{sv.23}
\end{eqnarray}

For this equation, we have that except the linear symmetry and the infinite
number of solution symmetries (we call them trivial symmetries) the equation admits
the vector field $Y_{2}=S\partial _{S}$, which is a reduced symmetry.
Therefore, the application of $Y_{2}$ to Equation (\ref{sv.23}) gives the second-order
ordinary differential equation
\begin{equation}
\beta ^{2}w_{,yy}+\left[ 2\alpha \left( m-y\right) +\frac{2\beta }{f\left(
y\right) }\left( \rho \kappa _{2}f^{2}\left( y\right) -\mu +r\right) \right]
w_{,y}+\left[ \left( \kappa ^{2}_{2}-1\right) f^{2}\left( y\right) +2\left(
r\kappa _{2}-r+\kappa _{1}\right) \right] w=0  \label{sv.24}
\end{equation}%
where $w=w\left( y\right) $ and
\begin{equation}
u\left( t,S,y\right) =S^{\kappa _{2}}\exp \left[ \kappa _{1}t\right] w\left(
y\right)   \label{sv.25}
\end{equation}%

Equation (\ref{sv.24}) is a linear second-order differential equation, and it
is well known that it is maximally symmetric and is invariant under the special linear (sl) algebra {$%
sl\left( 3,R\right) $} Lie algebra.

Similarly, if we perform a reduction with $Y_{2}$, the reduced equation
admits the Lie Symmetries $X_{1}$, $X_{v},~X_{b}$, and finally the solution
is again given by Equation (\ref{sv.25}) with the constraint Equation (\ref{sv.24})

Consider the application of the Lie symmetry vector $Y_{12}$ to Equation (\ref{sv.24}%
). We have that
\begin{equation}
u\left( t,S,y\right) =\exp \left[ k_{3}t\right] v\left( z,y\right)
~,~z=S\exp [-ct]  \label{sv.26}
\end{equation}%
where%
\begin{eqnarray}
0 &=&z^{2}f\left( y\right) v_{,zz}+2\rho \beta f\left( y\right)
v_{,zy}+\beta ^{2}v_{,yy}+2\left( r-c\right) zv_{,z}  \notag \\
&&+2\left[ \alpha \left( m-y\right) -\beta \rho \frac{\mu -r}{f\left(
y\right) }\right] v_{,y}-\left( r-\kappa _{3}\right) v  \label{sv.27}
\end{eqnarray}%

One can easily find that this equation only admits the Lie symmetry, $%
z\partial _{z}$, except the trivial symmetries, which is a reduced symmetry.
Therefore, the application of the zeroth-order invariants of the symmetry
vector $\left( z\partial _{z}+\kappa _{4}v\partial _{v}\right) $ in Equation (\ref%
{sv.27}) gives solution of the form Equation (\ref{sv.25}) with the constraint Equation (%
\ref{sv.24})

We continue with the determination of the group invariant solutions for
constant $f\left( y\right) $.

\subsection{Constant Volatility}

For $f\left( y\right) =f_{0}$, Equation (\ref{sv.10}) admits six Lie point
symmetries, plus the infinite number of solution symmetries. Moreover, Equation (\ref%
{sv.10}) is an $(1+2)$ evolution equation, and, in order to reduce it to an
ordinary differential equation, we have to apply the zeroth-order invariants
of two Lie symmetries. From Table \ref{table1}, we select reducing Equation (%
\ref{sv.10}) by using the following two-dimensional subalgebras $%
A_{I}=\left\{ Y_{1},Y_{2}\right\} $,~$B_{I}=\left\{ Y_{2},\bar{X}_{3}+\kappa
X_{u}\right\} $, \mbox{$B_{II}=\left\{ Y_{2},\bar{X}_{5}+\kappa X_{u}\right\} $}, \mbox{$%
C_{I}=\left\{ \bar{X}_{3}+\kappa X_{u},\bar{X}_{4}\right\} $}~and~\mbox{$%
C_{II}=\left\{ \bar{X}_{4},\bar{X}_{5}\right\} .$}

The reduction with the subalgebra $A_{I}$ we studied in the previous
subsection and the solution is Equation (\ref{sv.25}), where now from Equation (\ref{sv.24})
we have%
\begin{equation*}
w\left( y\right) =w_{1}\mathbf{K}\left( -\frac{c_{2}}{4\alpha },\frac{1}{2},%
\frac{\left( \alpha \left( m-y\right) +\frac{c_{1}}{2}\right) ^{2}}{\alpha }%
\right) +w_{2}\mathbf{\Lambda }\left( -\frac{c_{2}}{4\alpha },\frac{1}{2},%
\frac{\left( \alpha \left( m-y\right) +\frac{c_{1}}{2}\right) ^{2}}{\alpha }%
\right)
\end{equation*}%

We have that $w_{1}$ and $w_{2}$ are constants,
\begin{equation*}
c_{1}=\frac{2\beta }{f_{0}}\left( \rho \kappa _{2}f_{0}^{2}-\mu +r\right)
,~c_{2}=\left[ \left( \kappa ^{2}-1\right) f_{0}^{2}+2\left( r\kappa
_{2}-r+\kappa _{1}\right) \right]
\end{equation*}
and $\mathbf{K,\Lambda }$ are Kummer Functions.

We continue with the application of the remaining subalgebras.

The application of $B_{I}$ gives%
\begin{equation}
u\left( t,S,y\right) =S^{\kappa _{2}}\exp \left( \kappa e^{\alpha t}y\right)
\phi _{I}\left( t\right)
\end{equation}%
where $\phi _{I}\left( t\right) $ is given by the first-order ordinary
differential equation
\begin{equation}
2f_{0}\phi _{I,t}+\left[ \left( f_{0}^{3}\kappa _{2}+2f_{0}r\right) \left(
\kappa _{2}-1\right) +2\kappa \beta \left( f_{0}^{2}\rho \kappa _{2}-\mu +r+%
\frac{\alpha \mu f_{0}}{\beta }\right) e^{\alpha t}+f_{0}\beta ^{2}\kappa
^{2}e^{2\alpha t}\right] \phi _{I}=0
\end{equation}%
with solution%
\begin{equation}
\ln \left( \frac{\phi _{I}\left( t\right) }{\phi _{I0}}\right) =-\frac{%
\left( f_{0}^{2}\kappa _{2}+2r\right) }{2}\left( \kappa _{2}-1\right) t-%
\frac{2\kappa \beta }{2\alpha f_{0}}\left( f_{0}^{2}\rho \kappa _{2}-\mu +r+%
\frac{\alpha \mu f_{0}}{\beta }\right) e^{\alpha t}-\frac{\beta ^{2}\kappa
^{2}}{4\alpha }e^{2\alpha t}
\end{equation}

From the subalgebra $B_{II}$, we find the solution
\begin{equation}
u\left( t,S,y\right) =\phi _{II}\left( t\right) S^{\kappa _{2}}\exp \left(
\frac{\alpha }{\beta ^{2}}y^{2}+\left( \frac{\left( \kappa e^{-\alpha
t}-2\alpha m\right) }{\beta ^{2}f_{0}}-\frac{2\left( r-\mu +\rho \kappa
_{2}f_{0}^{2}\right) }{\beta f_{0}}\right) y\right)
\end{equation}%
where $\phi _{II}\left( t\right) $ is given by the expression
\begin{eqnarray}
\ln \left( \frac{\phi _{II}\left( t\right) }{\phi _{II0}}\right) &=&-\frac{1%
}{2}\left[ \left( f_{0}^{2}\kappa _{2}+2r\right) \left( \kappa _{2}-1\right)
+2\alpha \right] t \notag \\
&&+\left[ -\frac{\kappa }{\alpha \beta ^{2}f_{0}^{2}}\left( \beta \left(
r-\mu+\rho \kappa _{2}f_{0}^{2}\right) +\alpha mf_{0}\right) e^{-\alpha t}+%
\frac{\kappa ^{2}}{4\alpha f_{0}^{2}\beta ^{2}}e^{-2\alpha t}\right]
\end{eqnarray}

For the subalgebra $C_{I}$, we have the invariant solution
\begin{equation*}
u\left( t,S,y\right) =S^{\psi \left( t\right) }\exp \left( \kappa e^{\alpha
t}y\right) \phi _{III}\left( t\right)
\end{equation*}%
where
\begin{equation}
\psi \left( t\right) =-\frac{\kappa \beta \rho }{\alpha f_{0}t}-\frac{r}{%
f_{0}^{2}}+\frac{1}{2}+\frac{2}{4f_{0}^{2}t}
\end{equation}%
and $\phi _{III}\left( t\right) $ is given by the expression
\begin{eqnarray}
\ln \left( \frac{\phi _{III}\left( t\right) }{\phi _{III0}}\right) &=&-\frac{%
1}{2}\ln t+\frac{\left( 2\rho ^{2}-\alpha t\right) }{\alpha ^{2}t}\beta
^{2}\kappa ^{2}e^{2\alpha t}  \notag \\
&&-\frac{\kappa }{2\alpha f_{0}}e^{\alpha t}\left( 2\alpha mf+2\beta \left(
r-\mu -2\rho r+\frac{\rho }{2}f_{0}^{2}\right) +\frac{\left(
2r+f_{0}^{2}\right) ^{2}}{8f_{0}^{2}}t\right)
\end{eqnarray}

Finally, from the subalgebra $C_{II}$, we find the invariant solution%
\begin{equation}
u\left( t,S,y\right) =\frac{\phi _{IV}\left( t\right) }{\sqrt{at-2\rho ^{2}}}%
S^{\bar{\psi}\left( t,y\right) }\exp \left( V\left( t,y\right) \right)
\end{equation}%
where%
\begin{equation}
\bar{\psi}\left( t,y\right) =\frac{4\rho f_{0}a\left( m-y\right) +\beta
\left( \alpha \left( f_{0}^{2}-2r\right) t-4\rho \left( \mu -r\right)
\right) +\alpha \beta ^{2}}{2\beta f_{0}^{2}\left( \alpha t-2\rho
^{2}\right) }
\end{equation}%
\begin{equation}
V\left( t,y\right) =\left[ \frac{\alpha ^{2}}{\beta ^{2}\left( \alpha
t-2\rho ^{2}\right) }y-a\frac{2\beta \left( r-\mu +\frac{f_{0}^{2}}{2}\rho
-2r\rho \right) +2f_{0}\alpha m}{\beta ^{2}f_{0}^{2}\left( \alpha t-2\rho
^{2}\right) }\right] ty
\end{equation}
{and} function $\phi _{IV}\left( t\right) $ is given by the expression%
\begin{eqnarray}
\ln \left( \frac{\phi _{IV}\left( t\right) }{\phi _{IV0}}\right) &=&-\frac{%
f_{0}^{2}\left( 8\alpha -4r-f_{0}^{2}\right) -4r^{2}}{8f_{0}^{2}\left(
\alpha t-2\rho ^{2}\right) }\left( \alpha t^{2}-2\rho ^{2}t\right) +\frac{%
\left( f_{0}^{2}-2r\right) ^{2}}{2f_{0}^{2}\alpha \left( \alpha t-2\rho
^{2}\right) }\rho ^{4}  \notag \\
&&+\frac{2\left( f_{0}^{2}-2r^{2}\right) \left( f_{0}\alpha m-\beta \left(
\mu -r\right) \right) }{f_{0}^{2}\alpha \beta \left( \alpha t-2\rho
^{2}\right) }\rho ^{3}+\frac{\left( 4\beta \left( r-\mu \right)
+2f_{0}\alpha m\right) }{\beta ^{2}f_{0}\left( \alpha t-2\rho ^{2}\right) }%
m\rho ^{2}  \notag \\
&&+\frac{2\left( \mu -r\right) ^{2}}{f_{0}^{2}\alpha \left( \alpha t-2\rho
^{2}\right) }\rho ^{2}
\end{eqnarray}

In the following section, we study a special model for stochastic volatility
which has been proposed by Heston \cite{heston}.

\section{Heston Model}

\label{Heston}

In the Heston model for stochastic volatility the stock price, $S$, and the
volatility, $\sigma =f\left( y\right) $, satisfy the stochastic differential
equation given below%
\begin{eqnarray}
dS_{t} &=&rS_{t}dt+S_{t}\sqrt{Y_{t}}dW_{t}\quad \mbox{\rm and}  \label{hes.1},
\\
dY_{t} &=&\theta \left( \bar{m}-Y_{t}\right) dt+\delta \sqrt{Y_{t}}d\hat{Z}%
_{t}  \label{hes.2}
\end{eqnarray}%
where, in comparison with Equations (\ref{bs.00}) and (\ref{sv.03}), we observe that $%
f\left( y\right) =\sqrt{y}$ and $\beta =\delta \sqrt{y}$. The differential
equation which corresponds to that model is
\begin{eqnarray}
0 &=&\frac{1}{2}YS^{2}u_{,SS}+\rho \delta YSu_{,SY}+\frac{1}{2}\delta
^{2}Yu_{,YY}+rSu_{,S}  \notag \\
&&+\left( \theta \left( \bar{m}-Y\right) -\lambda Y\right) u_{,Y}-ru+u_{,t}
\label{hes.3}
\end{eqnarray}

Before we proceed with the symmetry analysis of Equation (\ref{hes.3}), we
perform the coordinate transformation $Y=y^{2}$. Then, Equation (\ref{hes.3})
becomes
\begin{eqnarray}
0 &=&\frac{1}{2}y^{2}S^{2}u_{,SS}+\beta \rho ySu_{,Sy}+\frac{1}{2}\beta
^{2}u_{,yy}+rSu_{,S}  \notag \\
&&+\frac{1}{2}\left[ c_{1}y+\left( \frac{c_{2}}{y}\right) \right]
u_{,y}-ru+u_{,t}  \label{hes.4}
\end{eqnarray}%
in which we have made the replacements $\frac{\delta }{2}\rightarrow \beta
,~c_{1}=\left( \theta -\lambda \right) $ and $c_{2}=\theta m-\beta ^{2}$. \
Equation (\ref{hes.4}) can be compared with Equation (\ref{sv.05}) for $f\left(
y\right) =y$. However, the risk premium factor of Equation (\ref{sv.09}) is not zero
and has absorbed the term of the Orstein--Uhlenbeck process.

From the Lie symmetry condition, Equations (\ref{sv.17}) for (\ref{hes.4}) we
find that this equation admits the \mbox{Lie symmetries}
\begin{equation}
X_{1}=\partial _{t}~,~X_{2}=S\partial _{S}~  \label{hes.5}
\end{equation}%
plus the $X_{u},~X_{b}$, that is, Equation (\ref{hes.4}) is invariant under
the Lie algebra $\left\{ 3A_{1}\right\} \oplus _{s}\infty A_{1}$.

We continue with the application of the Lie symmetries in order to reduce
Equation (\ref{hes.5}). \mbox{We follow} the results of Section \ref{arbi}, that
is, we apply the group invariants of the subalgebra $\left\{
Y_{1},Y_{2}\right\} $.

We find that the corresponding invariant solution of Equation (\ref{hes.4}) is
\begin{equation}
u\left( t,S,y\right) =S^{\kappa _{2}}e^{\kappa _{1}t}W\left( y\right)
\label{hes.6}
\end{equation}%
where $W\left( y\right) $ satisfies the linear second-order differential
equation:%
\begin{equation}
\beta ^{2}W_{,yy}+\left( 2\kappa _{2}\beta \rho +c_{1}y+\frac{c_{2}}{y}%
\right) W_{,y}+\left( 2\kappa _{1}-2r\left( 1-\kappa _{2}\right)
+y^{2}\left( \kappa _{2}^{2}-\kappa _{1}^{2}\right) \right) W=0
\label{hes.7}
\end{equation}%
the solution in closed form of which is expressed in terms of Kummer
Functions.

\begin{figure}[tpb]
\centering
\includegraphics[width=12cm,height=8cm]{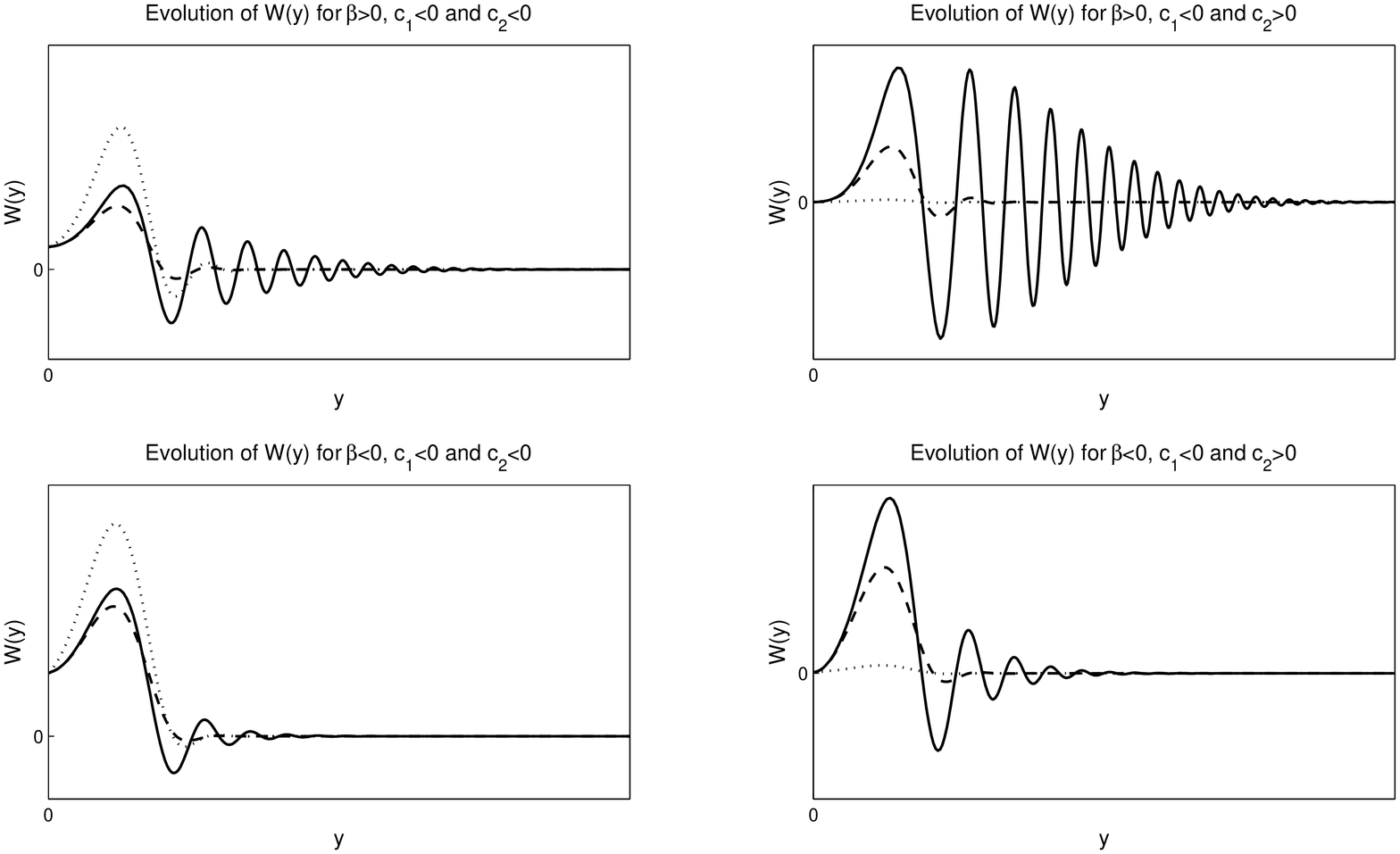}
\caption{Evolution of the solution Equation (\protect\ref{hes.7}) of the Heston
model. For the numerical solutions, we select $\protect\rho =0.5,~\protect%
\beta =\left\vert 0.7\right\vert ,r=0.5,~\protect\kappa _{1}=1~$and $\protect%
\kappa _{2}=0.5$. The left figures are for negative $c_{1},~c_{2}$, while
the right figures are for negative $c_{1}$ and positive $c_{2}$. The solid
lines are for $\left\vert c_{2}\right\vert =5\left\vert c_{1}\right\vert ,~$%
the dotted lines are for $\left\vert c_{2}\right\vert =0.2\left\vert
c_{1}\right\vert ,$ and the dash-dash lines are for $\left\vert
c_{2}\right\vert =\left\vert c_{1}\right\vert $. \ The top figures are for $%
\protect\beta >0$, while the lower figures for $\protect\beta <0$.}
\label{plot1}
\end{figure}

\begin{figure}[tpb]
\centering
\includegraphics[width=12cm,height=8cm]{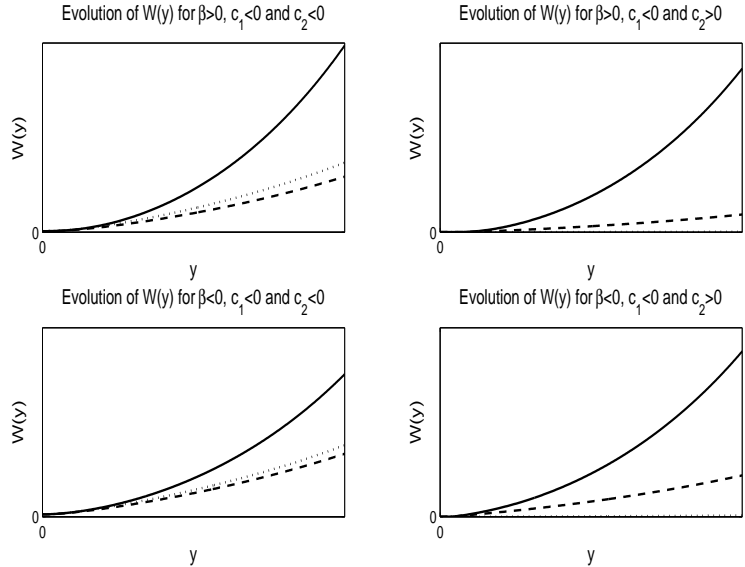}
\caption{Evolution of the solution Equation (\protect\ref{hes.7}) of the Heston
model. For the numerical solutions, we select $\protect\rho =0.5,~\protect%
\beta =\left\vert 0.7\right\vert, r=0.5,~\protect\kappa _{1}=1~$and $\protect%
\kappa _{2}=1.01$. The left figures are for negative $c_{1},~c_{2}$, while
the right figures are for negative $c_{1}$ and positive $c_{2}$. The solid
lines are for $\left\vert c_{2}\right\vert =\frac{3}{2}\left\vert
c_{1}\right\vert ,~$the dotted lines are for $\left\vert c_{2}\right\vert =%
\frac{2}{3}\left\vert c_{1}\right\vert ,$ and the dash-dash lines are for $%
\left\vert c_{2}\right\vert =\left\vert c_{1}\right\vert $. \ The top
figures are for $\protect\beta >0$, while the lower figures for $\protect%
\beta <0$.}
\label{plot3}
\end{figure}

In Figures \ref{plot1} and \ref{plot3}, we give numerical solutions of
Equation (\ref{hes.7}). Figure \ref{plot1} is for negative value of $\Delta
\kappa ,~$whereas Figure \ref{plot3} ~is for positive value of $\Delta
\kappa $, where $\Delta \kappa =\kappa _{2}-\kappa _{1}$.

\vspace{-12pt}

\section{Stein--Stein Model}

\label{steinstein}

The model which has been proposed by Elias M. Stein\ and Jeremy C. Stein \cite{steinstein}
\ describes an European option with stochastic volatility for which the
correlation among the two Brownian motions vanishes, \textit{i.e.}, \mbox{$\rho =0$},
in Equation (\ref{sv.04}). Moreover, they considered that the risk premium factor is
constant, \textit{i.e.}, $\gamma \left( t,S,y\right) =\gamma _{0}$ and the
volatility is $f\left( y\right) =y$, while the stochastic differential
equation is%
\begin{eqnarray}
dS_{t} &=&rS_{t}dt+S_{t}Y_{t}dW_{t}\quad \mbox{\rm and} \\
dY_{t} &=&\theta \left( \bar{m}-Y_{t}\right) dt+\delta d\hat{Z}_{t}
\end{eqnarray}

Therefore, from Equation (\ref{sv.05}), we have that the $(1+2)$ evolution differential
equation of the Stein--Stein model is
\begin{equation}
\frac{1}{2}y^{2}S^{2}u_{,SS}+\frac{1}{2}\beta ^{2}u_{,yy}+rSu_{,S}+\left[
\alpha \left( m-y\right) -\beta \gamma _{0}\right] v_{,y}-ru+u_{,t}=0
\label{stein.01}
\end{equation}

From the Lie symmetry condition Equation (\ref{sv.17}), we find that Equation (\ref%
{stein.01}) admits the \mbox{Lie symmetries}%
\begin{equation}
X_{1}=\partial _{t},~X_{2}=S\partial _{S}~,~X_{u}=u\partial
_{u},~X_{b}=b\left( t,S,y\right) \partial _{u}  \label{stein.02}
\end{equation}%
which form the Lie algebra $\left\{ 3A_{1}\right\} \oplus _{s}\infty A_{1}$,
and it is the admitted algebra of the Heston model and of Equation (\ref%
{sv.10}) for the arbitrary function $f\left( y\right) $.

Following the steps of the previous sections, we find that the invariant
solution of the Stein--Stein model with respect to the Lie algebra $\left\{
Y_{1},Y_{2}\right\} $ is
\begin{equation}
u\left( t,S,y\right) =S^{\kappa _{2}}e^{\kappa _{1}t}Y\left( y\right)
\label{stein.03}
\end{equation}%
where $Y\left( y\right) $ is given by the linear second-order differential
equation
\begin{equation}
\beta ^{2}Y_{,yy}+\left( 2\omega -2ay\right) Y_{,y}+\left( 2r\left( \kappa
_{2}-1\right) +2\kappa _{1}+y^{2}\left( \kappa _{2}^{2}-\kappa
_{1}^{2}\right) \right) Y=0  \label{stein.04}
\end{equation}

\begin{figure}[tpb]
\centering
\includegraphics[width=12cm,height=8cm]{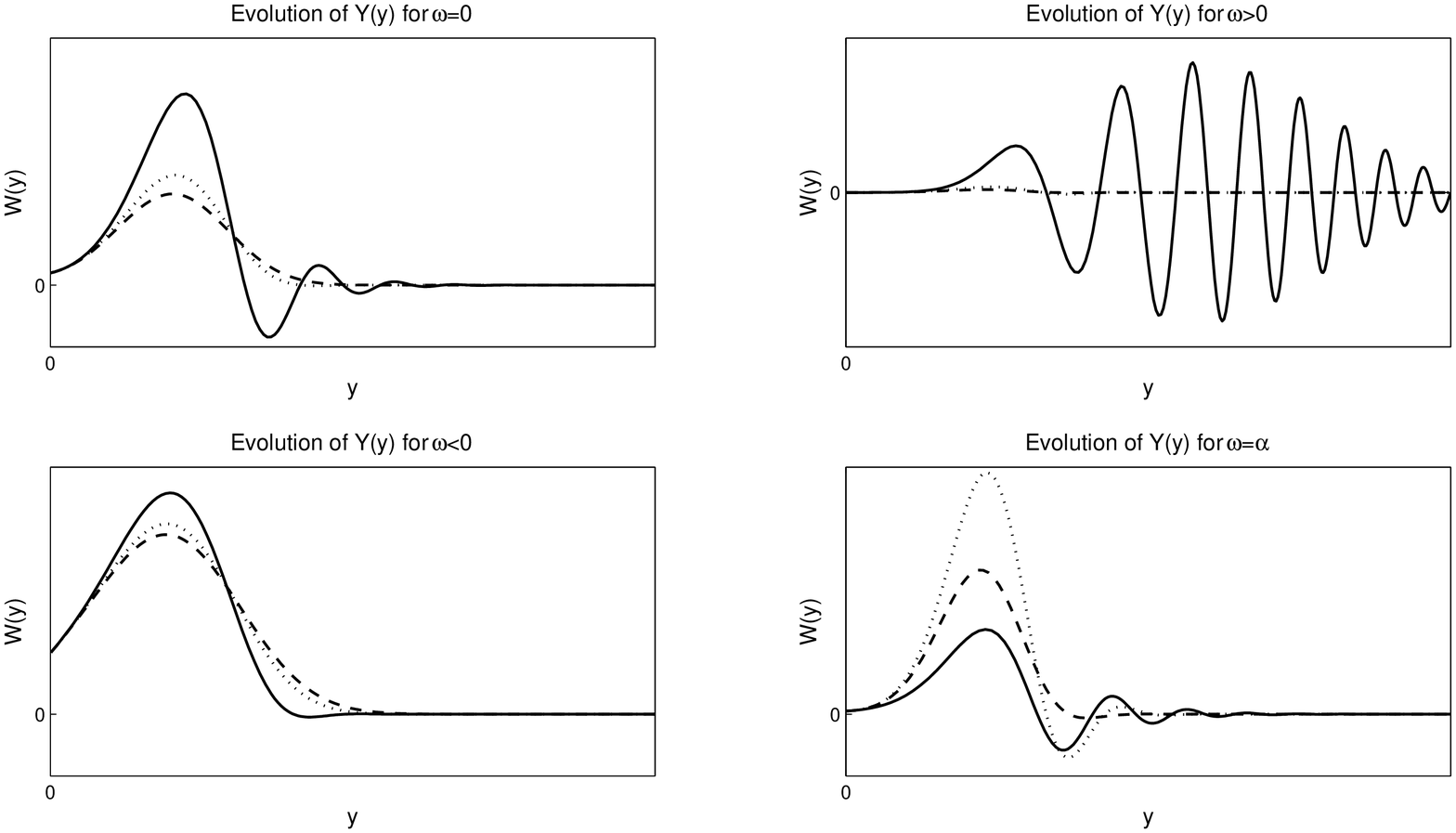}
\caption{Evolution of the solution Equation (\protect\ref{stein.04}) of the
Stein--Stein model. For the numerical solutions, we select $~\protect\beta %
=0.5,~r=0.5,~\protect\kappa _{1}=1~$and $\protect\kappa _{2}=0.5$. The
figures are for $\protect\omega =0,~\protect\omega =0.5,~\protect\omega %
=-0.5 $ and $\protect\omega =\protect\alpha, ~$respectively. The solid lines
are for $\protect\alpha =0.1,~$the dotted lines are for $\protect\alpha =0.3$
and the dash-dash lines are for $\protect\alpha =0.4$. }
\label{stein1}
\end{figure}

\begin{figure}[tpb]
\centering
\includegraphics[width=15cm,height=9cm]{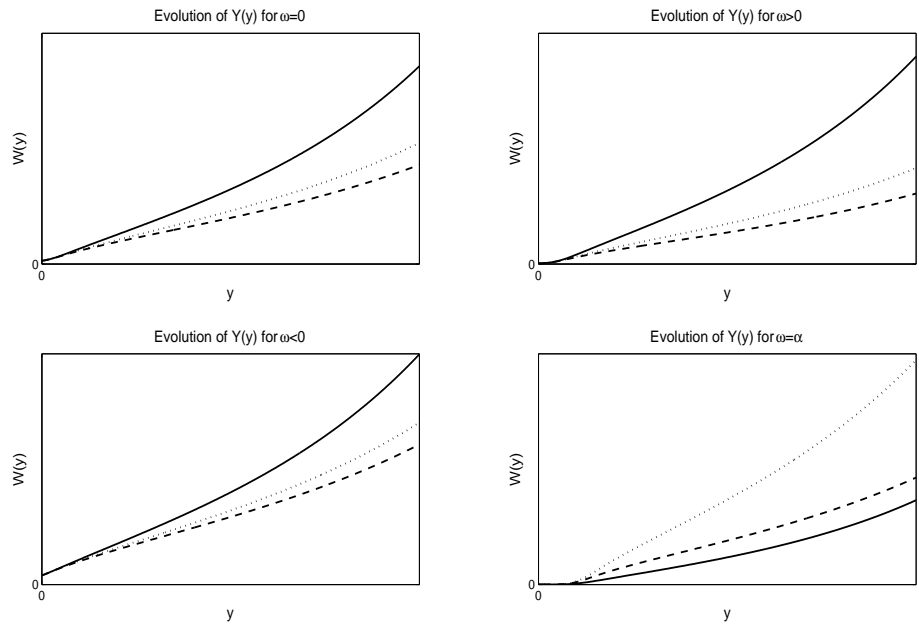}
\caption{Evolution of the solution Equation (\protect\ref{stein.04}) of the
Stein--Stein model. For the numerical solutions, we select $~\protect\beta %
=0.5,~r=0.5,~\protect\kappa _{1}=1~$and $\protect\kappa _{2}=1.01$. The
figures are for $\protect\omega =0,~\protect\omega =0.5,~\protect\omega %
=-0.5 $ and $\protect\omega =\protect\alpha $, respectively. The solid lines
are for $\protect\alpha =1.1,~$the dotted lines are for $\protect\alpha =1.3$
and the dash-dash lines are for $\protect\alpha =1.4$. }
\label{stein2}
\end{figure}

The closed-form solution of this equation can be expressed in terms of the
Hypergeometric Functions, where $\omega =\alpha m-\beta \gamma _{0}$. In
Figures \ref{stein1} and \ref{stein2}, we give the numerical evolution of $%
Y\left( y\right) $ for {various values of the parameters,} $\omega ~$and $\alpha $,
for negative $\Delta \kappa ~$and positive $\Delta \kappa, ~$%
respectively, where $\Delta \kappa =\kappa _{2}-\kappa _{1}$.
\vspace{-12pt}

\section{Conclusions}

\label{conclusions}

Volatility with a stochastic process has been shown to be essential for
Financial Mathematics. In this work, we studied the algebraic properties,
\textit{i.e.}, the Lie symmetries, of the modified Black--Scholes--Merton
Equation for European options with a stochastic volatility. We have shown
that the autonomous model without the risk premium factor is invariant under
a group of point transformations which form the $\left\{ 3A_{1}\right\}
\oplus _{s}\infty A_{1}$ Lie Algebra for an arbitrary functional form of the
volatility, $\sigma $. Moreover, when the volatility is constant but the
price of the option depends on the second Brownian motion, in which the
volatility is defined, the modified Black--Scholes--Merton Model is invariant
under six, plus the infinity, Lie point symmetries and it is not maximally
symmetric as the$~$Black--Scholes--Merton Equation with nonstochastic
volatility is.

Furthermore, we showed that the Black--Scholes--Merton Model, in which the
volatility is constant, but is defined by an Orstein--Uhlenbeck
process, is invariant under the same group of point transformations as
that of the two-factor model of commodities. The reason for that is that the
two models have in common the terms which follow from the
Orstein--Uhlenbeck process.

Moreover, we applied the zeroth-order invariants of the Lie symmetries, and we
reduced the model to a linear second-order differential equation. As far as
the case of constant volatility is concerned, we found the closed forms of
the group invariant solutions.

Finally, we studied the algebraic properties and the invariant solutions of
two models, the Heston model and the Stein--Stein model, with stochastic
volatility of special interest. For each model, we found the invariant
solution and we gave some figures for the evolution of the models. Of course
because Equation (\ref{sv.05}) is a linear equation, the general solution is given by
the linear combination of the invariant solutions that we have found, while
the latter are constrained by the initial conditions and the boundary
conditions of the model.

A general consideration of Equation (\ref{sv.05}), in which the risk premium
factor plays a role is still in progress, and the results will be
published elsewhere.


\vspace{6pt}
\textbf{Aacknowledgments:} The research of Andronikos Paliathanasis was supported by FONDECYT grant No. 3160121. K. Krishnakumar thanks the University
Grants Commission for providing a University Grants Commission-Basic Scientific Research Fellowship to perform this research work.

\bibliographystyle{mdpi}




\end{document}